\pdfoutput=1
\RequirePackage{ifpdf}
\ifpdf 
\documentclass[pdftex]{sigma}
\else
\documentclass{sigma}
\fi

\usepackage[all]{xy}

\newtheorem{Theorem}{Theorem}[section]
\newtheorem{Corollary}[Theorem]{Corollary}
\newtheorem{Lemma}[Theorem]{Lemma}
\newtheorem{Proposition}[Theorem]{Proposition}
{\theoremstyle{definition}
\newtheorem{Definition}[Theorem]{Definition}

\newtheorem{Example}[Theorem]{Example}
\newtheorem{Remark}[Theorem]{Remark}
}

\newcommand{\kth}{k^{\mathrm{th}}}

\newcommand{\dd}{\mathrm{d}}

\newcommand{\R}{\mathbb{R}}

\newcommand{\JkXM}{J^k(X, M)}

\newcommand{\ttt}{\mathbf{t}}
\newcommand{\ttv}{\mathbf{v}}
\newcommand{\nn}{\mathbf{n}}
\newcommand{\Act}{\mathcal{A}}

\begin{document}

\allowdisplaybreaks

\renewcommand{\thefootnote}{$\star$}

\renewcommand{\PaperNumber}{036}

\FirstPageHeading

\ShortArticleName{On Local Congruence of Immersions in Homogeneous or Nonhomogeneous Spaces}

\ArticleName{On Local Congruence of Immersions\\ in Homogeneous or Nonhomogeneous Spaces\footnote{This paper is a~contribution to the Special Issue  ``Symmetries of Dif\/ferential
Equations:  Frames, Invariants and~Applications''.
The full collection is available at
\href{http://www.emis.de/journals/SIGMA/SDE2012.html}{http://www.emis.de/journals/SIGMA/SDE2012.html}}}

\Author{Jeongoo CHEH}

\AuthorNameForHeading{J.~Cheh}
\Address{Department of Mathematics \& Statistics, The University of Toledo, Toledo, OH 43606, USA}

\Email{\href{mailto:jeongoocheh@gmail.com}{jeongoocheh@gmail.com}}

\ArticleDates{Received May 14, 2012, in f\/inal form April 19, 2013; Published online April 28, 2013}

\Abstract{We show how to f\/ind a complete set of necessary and suf\/f\/icient conditions that solve the f\/ixed-parameter local congruence problem of immersions in $G$-spaces, whether homogeneous or not, provided that a certain $\kth$ order jet bundle over the $G$-space admits a $G$-invariant local coframe f\/ield of constant structure.  As a corollary, we note that the dif\/ferential order of a minimal complete set of congruence invariants is bounded by $k+1$.  We demonstrate the method by rediscovering the speed and curvature invariants of Euclidean planar curves, the Schwarzian derivative of holomorphic immersions in the complex projective line, and equivalents of the f\/irst and second fundamental forms of surfaces in $\R^3$ subject to rotations.}

\Classification{53A55; 53B25}

\begin{flushright}
\begin{minipage}{80mm}
\it  Dedicated to my teacher, Professor Peter Olver,\\
in honor of his sixtieth birthday.
\end{minipage}
\end{flushright}

\renewcommand{\thefootnote}{\arabic{footnote}}
\setcounter{footnote}{0}

\section{Introduction}
The equivalence problem of immersed submanifolds of a manifold $M$ is to f\/ind a set of computable criteria that determine whether the submanifolds, or open subsets thereof in their own topology, are equivalent or not under the action of a prescribed symmetry group\footnote{The {\em symmetry group} of a space is by def\/inition the group of all admissible transformations that are typically characterized to preserve, for instance, certain geometric structures on the space; conversely, prescribing a symmetry group to a space amounts to declaring the class of admissible transformations of the space.}~$G$ of~$M$.  Recall that a~manifold~$M$ on which a (local) Lie group $G$ acts is called a $G$-space.  For a homogeneous $G$-space~$M$, there is a classical solution to the equivalence problem, known as Cartan's moving frame method,~\cite{green, griffiths, landsberg, sternberg}, that consists of lifting the submanifolds into the principal $H$-bundle~$G$ over~$M$, where $H < G$ is the stabilizer of a point of~$M$, to pull the Maurer--Cartan forms of~$G$ back down to the submanifolds and then use an argument resembling {\em Cartan's technique of the graph}, \cite{warner}, that is Cartan's rendition of the Frobenius theorem.

To our knowledge, however, little has been known of any general solution for the case of submanifolds of {\em nonhomogeneous} spaces.  One of the obvious obstacles in one's attempt to apply the classical method to the case of a nonhomogeneous space $M$ is the unavailability of any natural manifestation of $G$ as a principal bundle over $M$, unlike in the case of homogeneous spaces.

In \cite{fels-olver1, fels-olver2}, remarks (for example, Theorem~7.2 in \cite{fels-olver2})
have been made on a solution to a certain type of equivalence problems of submanifolds of a space where its symmetry group may not act transitively.  Its method involves using an {\em extended invariant coframe field} that consists of an invariant coframe f\/ield along with a complete system of dif\/ferential invariants.
Although the papers dealt with equivalence problems somewhat dif\/ferent\footnote{Shortly, we will  explain two dif\/ferent types of equivalence problems: the {\em variable-parameter} and {\em fixed-parameter} types.} from those that we will study in our own paper, they have stimulated our curiosity and motivation by drawing our attention to the article~\cite{griffiths} that studied equivalence problems for submanifolds immersed exclusively in homogeneous spaces.

Besides the issue of whether the ambient space is homogeneous or nonhomogeneous, there is another technicality related to {\em reparametrization} of submanifolds.  Suppose that \mbox{$\psi_1, \psi_2: X{\longrightarrow} M$} are two submanifolds of a $G$-space~$M$.  One type of equivalence problem is to determine whether there exits an element $g\in G$ that transforms one image set $\psi_1(X)$ to the other image set~$\psi_2(X)$, i.e., $g(\psi_1(X))=\psi_2(X)$.  Since this problem is concerned only with the images of the submanifold immersions, reparametrizing a submanifold, such as $\psi_1\circ \phi:X\longrightarrow M$ with a (local) dif\/feomorphism $\phi:X\longrightarrow X$, should not change the f\/inal answer of the equivalence problem.

The equivalence problem allowing reparametrization of submanifolds has been studied extensively, for instance, in~\cite{equivalence-book}.  The overall idea of the solution for this problem largely consists of two parts. First, f\/ind a (minimal) set of dif\/ferential invariants that generate, through invariant dif\/ferentiation, the entire algebra of dif\/ferential invariants.  Then restrict the dif\/ferential invariants of up to a certain enough order to the given submanifolds to obtain {\em signature submanifolds}, which in turn leads to the f\/inal question on {\em overlapping submanifolds}.  We should note the theoretical development in a novel and systematic method of f\/inding a generating set of dif\/ferential invariants under f\/inite- or inf\/inite-dimensional Lie pseudo-group actions, \cite{dia, opmf}, and its applications to a number of dif\/ferent geometries, \cite{hubert, surfaces, centro-affine}.

The other type of equivalence problem is concerned with not only the images of the immersions but also the immersions themselves as well.  In other words, even when the images may be regarded equivalent in the previous sense (that we call {\em variable-parameter congruence}) of allowing reparametrization of immersed submanifolds, if one of the submanifolds is indeed nontrivially reparametrized, the immersed submanifolds are no longer to be considered equivalent in the new sense.  It is this type of more restrictive (and hence easier) equivalence problem that we will study in our paper, and we call it the {\em fixed-parameter congruence problem}.
The precise meaning of congruence is as follows:
\begin{Definition}\label{define congruence}
Suppose that $\psi_1, \psi_2: X\longrightarrow M$ are immersions of a manifold $X$ into a $G$-space~$M$.  We say that $\psi_1$ and $\psi_2$ are {\em congruent at $x_0\in X$} if there exist an open neighborhood $U\subset X$ of $x_0$ and a transformation $g\in G$ such that $\psi_1(x)=g\circ \psi_2(x)$ for all $x\in U$.
\end{Definition}
The main goal of our paper is to provide a theoretical justif\/ication of a method that we devised for solving the local congruence problem of immersed submanifolds of a $G$-space that may or may not be homogeneous.  In essence, our method is a hybrid of the classical Cartan's method, \cite{green, griffiths,landsberg, sternberg}, and the method of invariant coframe f\/ields constructed by equivariant moving frames, \cite{fels-olver1, fels-olver2, kogan-olver, opmf}.  On the one hand, we follow the classical idea of utilizing Cartan's technique of the graph which provides a means to solve the equivalence problem by examining certain invariant coframe f\/ields;  on the other hand, we use the machinery of equivariant moving frames, as developed by Olver et al., to construct particular invariant coframe f\/ields that are of {\em constant structure} (Def\/inition~\ref{define constant-structure coframes}).  One of our main results is Theorem~\ref{general solution} where we give a~proof to the ef\/fect that we can use an invariant coframe f\/ield of constant structure, instead of the unavailable Maurer--Cartan forms of homogeneous spaces required by the classical method, to completely determine congruence of immersions whether the ambient space is homogeneous or not.
We regard our result as extending and generalizing to quite arbitrary $G$-spaces the key lemmas proved in \cite{griffiths} that solved congruence problems in homogeneous spaces.  A notable corollary to Theorem~\ref{general solution} is that if the prolongation of the ({\em verticalized}\footnote{See Def\/inition~\ref{define verticalized actions}.}) action of the group becomes locally free at order~$k$, then a minimal complete set of congruence invariants should be of order $k+1$ or less.

We illustrate our method by f\/inding congruence conditions of some specif\/ic well-known examp\-les, which will thus help attest to the validity of our method.  Through the f\/irst few sections of the paper, we explain key ideas of our method with the example of the congruence problem of curves in $\R^2$ under the action of the orientation-preserving rigid motion group ${\rm SE}(2)$, whereby we obtain speed and curvature congruence invariants of curves.  The last section of the paper focuses on three examples.  We rediscover the  Schwarzian derivative as a congruence invariant of holomorphic immersions, or biholomorphisms, in the complex projective line $\mathbb{CP}^1$ under the action of the projective special linear group ${\rm PSL}(2,\mathbb{C})$ or alternatively in the complex plane~$\mathbb{C}$ under the action of linear fractional, or M\"obius,  transformations.  The next example obtains a complete system of congruence invariants, {\em including} equivalents to the f\/irst and second fundamental forms of surfaces in~$\R^3$, under the intransitive action of the rotation group ${\rm SO}(3)$.  The last example computes congruence conditions for surfaces in $\R^3$ under another intransitive action of a certain subgroup of a Heisenberg group.

Throughout the paper, all manifolds and maps are assumed to be smooth, $M$ always denotes a smooth manifold, and $G$ an $r$-dimensional local Lie group acting on~$M$ and other manifolds.  Furthermore, since we are concerned only with {\em local} equivalence, by a manifold will we actually mean a connected open subset thereof that is small enough to suit the context, and likewise by~$G$ a suf\/f\/iciently small connected open neighborhood of its identity element.  However, whenever serious confusion may be possible, we will explicitly remind ourselves of the local nature of various spaces and maps involved.

\section{Verticalization and prolongation of group actions}

To study the ef\/fects of group actions on immersed submanifolds of $M$, we begin by treating the submanifolds as (local) sections of the product bundle $X\times M\longrightarrow X$ where $X$ denotes the space of parameters of the submanifolds or the domain of the maps immersing the submanifolds.  This viewpoint establishes a one-to-one correspondence $\gamma$ between the set of immersions $f:X\longrightarrow M$ and the set of sections $\gamma(f)$ def\/ined by
\[
\gamma(f): \ X\longrightarrow X\times M, \qquad \gamma(f)(x)=(x,f(x)).
\]

The action of $G$ on $M$ induces an obvious natural action on $X\times M$:
\begin{Definition}\label{define verticalized actions}
The {\em verticalized action} of $G$ is def\/ined to be
\[
G\times (X\times M)\longrightarrow X\times M, \qquad (g,(x,u))\longmapsto (x,g\cdot u).
\]
\end{Definition}

It is this action of $G$ on $X\times M$ that we will later prolong to the jet bundles $\JkXM$, $k=0,1,2,\dots$,  of local sections of the bundle $X\times M\longrightarrow X$.

\begin{Remark}
In the literature, if $M$ itself is a f\/ibered manifold, it is customary, for the sake of computation in local coordinates, to consider the jet bundles $J^kM$ of only those submanifolds of $M$ that are transverse to the f\/ibers of $M$.  Our construction of the bundles $\JkXM$, on the other hand, contains prolongations, \cite{app-of-Lie-groups, equivalence-book}, of {\em all} submanifolds of $M$, including even those that may not be transverse to the f\/ibers of $M$.  Furthermore, the local coordinate description of the prolonged action of $G$ on $\JkXM$ turns out much simpler than that on $J^kM$ since the verticalized action of $G$ on $X\times M$ does not af\/fect base ($X$) coordinates at all.
\end{Remark}

Let $(x^1, \dots, x^p)$ and $(u^1,\dots ,u^n)$ denote local coordinate systems on manifolds $X$  and $M$, respectively.  We write $(x^1, \dots, x^p, u^1,\dots, u^n, \dots, u^\alpha_J, \dots )$, where $\alpha=1,2,\dots, n$ and $J$ is a symmetric multi-index over $\{1,2,\dots, p\}$, for the standard bundle-adapted local coordinate system on $\JkXM$.  Then, for $g\in G$ and $z\in \JkXM$, while the base coordinates of $g\cdot z$ remain unchanged:
\begin{displaymath}
x^i(g\cdot z)=x^i(z), \qquad i=1,2,\dots, p,
\end{displaymath}
all the f\/iber coordinates $u^\alpha_{J}(g\cdot z)$ are found by iterating the recursive formula
\begin{gather}\label{vertical coord}
u^\alpha_{J,i}(g\cdot z)=D_i \big(u^\alpha_{J}(g\cdot z)\big), \qquad \alpha=1,2,\dots , n, \qquad i=1,2,\dots, p,
\end{gather}
where
\[
D_i:=\frac{\partial}{\partial x^i}+\sum_{|J|\geq 0}\sum_{\alpha=1}^n u^\alpha_{J,i}\frac{\partial}{\partial u^\alpha_J}, \qquad i=1,2,\dots, p,
\]
are the coordinate-wise {\em total differential operators} on the inf\/inite-order Jet bundle $J^\infty(X, M)$.  For more details on calculus on jet bundles and related notations, refer, for example, to the resources \cite{anderson, app-of-Lie-groups, equivalence-book}.

\begin{Example}[Euclidean plane]\label{SE(2) action original}
Consider the Euclidean plane $(\R^2, {\rm SE}(2))$ under the transitive action of the group  ${\rm SE}(2):={\rm SO}(2)\ltimes \R^2$ of orientation-preserving rigid motions.  Let $\theta\in\R$ and $(a,b)\in\R^2$ be the parameters of ${\rm SE}(2)$ such that, for $(u,v)\in\R^2$, the action is given by
\[
(\theta,a,b)\cdot (u,v)\longmapsto (\widehat{u},\widehat{v}),
\]
where

\[
\widehat{u}=u\cos\theta-v\sin\theta+a, \qquad \widehat{v}=u\sin\theta+v\cos\theta+b.
\]
Let $x$ and $(u,v)$ denote the local coordinate systems on $\R$ and $\R^2$, respectively.   Then the jet bundle $J^1(\R,\R^2)$ of local sections $\gamma(\psi)$ of the product bundle $\R\times\R^2\longrightarrow \R$,
where $\psi:\R\longrightarrow \R^2$ is a locally def\/ined immersion,
has the standard local coordinate system $(x,u,v,u_x,v_x)$.  In these coordinates, the verticalized action of ${\rm SE}(2)$ on $\R\times\R^2=J^0(\R,\R^2)$ is given by the rule
\[
(\theta,a,b)\cdot (x,u,v)\longmapsto (x,\widehat{u},\widehat{v})
\]
and its prolonged action $\mathcal{A}:{\rm SE}(2)\times J^1(\R,\R^2)\longrightarrow J^1(\R,\R^2)$  by
\[
(\theta,a,b)\cdot (x,u,v,u_x,v_x)\longmapsto (x,\widehat{u},\widehat{v},\widehat{u_x}, \widehat{v_x})
\]
where, according to (\ref{vertical coord}),
\[
\widehat{u_x}=D_x\widehat{u}=u_x\cos\theta-v_x\sin\theta, \qquad \widehat{v_x}=D_x\widehat{v}=u_x\sin\theta+v_x\cos\theta.
\]
\end{Example}

\section{Review of equivariant moving frames}
The purpose of this section is to review and use equivariant moving frames as part of the tools that we will need later on to solve congruence problems.  Our main references are \cite{fels-olver1,fels-olver2}, and we will provide proofs for all our claims whenever they do not have exact counterparts proved in the references.

In general, two special types of maps between $G$-spaces stand out in the business of constructing invariant dif\/ferential forms.
Let $\psi:M\longrightarrow N$ be a locally def\/ined map between two $G$-spaces $M$ and $N$.
\begin{Definition}\label{invariance-equivariance}
Let $V\subset G$ be an open neighborhood of the identity element of $G$.  The map $\psi$ is said to be {\em $G$-invariant} if
\[
g^*\psi:=\psi\circ g=\psi
\]
for all $g\in V$ such that $\psi\circ g$ is well-def\/ined; and {\em $G$-equivariant} if
\begin{gather}\label{define equivariance}
\psi(g\cdot x)=g\cdot\psi(x)
\end{gather}
for all $g\in V$ and $x\in M$ such that the both sides of the equation (\ref{define equivariance}) are well def\/ined.
\end{Definition}

\begin{Remark}
Throughout the paper, all our constructions, without exception, of maps, dif\/fe\-ren\-tial forms, group actions, etc. will be done locally. To avoid the trite use of the word {\em local} and relevant open subset notations, we henceforth impose the blanket assumption that all domains of maps should be understood to be freely replaceable by appropriately small open subsets thereof in order to make sense of the expressions in which they are considered.  Thus, for example, we will simply say that a map $\psi$ is $G$-equivariant instead of more precisely saying that $\psi$ is $V$-equivariant for an open neighborhood $V\subset G$ of the identity element of $G$, and this is what is ref\/lected in Def\/inition~\ref{invariance-equivariance}.
\end{Remark}

\begin{Proposition}\label{equivariant-invariant}
If $\psi$ is $G$-invariant and $\omega$ is any differential form on $N$, then $\psi^*\omega$ is a~$G$-invariant differential form on $M$.  On the other hand, if $\psi$ is $G$-equivariant and $\omega$ is a $G$-invariant differential form on $N$, then $\psi^*\omega$ is a $G$-invariant differential form on $M$.
\end{Proposition}
\begin{proof}
If $\psi:M\longrightarrow N$ is $G$-invariant, then for any $g\in G$ and any dif\/ferential form $\omega$ on $N$, $g^*(\psi^*\omega)=(\psi\circ g)^*\omega=\psi^*\omega=\omega$, showing that $\psi^*\omega$ is $G$-invariant.  The proof for the case of equivariant maps is similar.
\end{proof}
Thus, when $N$ is the Lie group $G$ that always comes with a canonical family of invariant dif\/ferential forms, (that is, Maurer--Cartan forms), a $G$-equivariant map $M\longrightarrow G$ plays a key role in producing invariant dif\/ferential forms on $M$, and thus has acquired a special recognition in~\cite{fels-olver2}.
\begin{Definition}\label{define moving frames}
Under an action of $G$ on itself, a locally def\/ined $G$-equivariant map $\rho:M{\longrightarrow}G$, is called an {\em $($equivariant$)$ moving frame}.  Associated to $\rho$ is the locally def\/ined {\em $($equivariant$)$ moving frame section}
\[
\sigma: \ M\longrightarrow G\times M, \qquad z\longmapsto (\rho(z),z),
\]
of the trivial product bundle $G\times M\longrightarrow M$.
\end{Definition}

When the $G$-action on $M$ is locally free, a typical procedure for obtaining a moving frame begins by choosing a local cross-section $\Gamma$ transverse to the $G$-orbits in $M$ and then def\/ining a~map $\rho:M\longrightarrow G$ by requiring $\rho(z)\cdot z\in \Gamma$ for $z\in M$.  If the $G$-action on $M$ is a left action\footnote{An action $G\times M\longrightarrow M$, $(g,z)\longmapsto g\cdot z$, is called a {\em left} action if $g_1\cdot (g_2\cdot z)=(g_1g_2)\cdot z$ for $g_1, g_2\in G$ and $z\in M$; or a {\em right} action if $g_1\cdot (g_2\cdot z)=(g_2g_1)\cdot z$, where $g_1g_2$ and $g_2g_1$ simply signify the multiplication structure of the group~$G$.} and if $G$ acts on itself by the rule $g\cdot h\longmapsto hg^{-1}$, then the map $\rho$ def\/ined in such a way is indeed a~moving frame satisfying the condition: for all $z\in M$ and $g\in G$, $\rho(g\cdot z)=\rho(z) g^{-1}$ whenever the two sides are well def\/ined. If we have a right $G$-action on~$M$ and if we def\/ine the right action of $G$ on itself by $g\cdot h\longmapsto g^{-1}h$, then the map $\rho$ constructed above is again a moving frame, this time satisfying $\rho(g\cdot z)=g^{-1}\rho(z)$.

As long as $G$ acts on itself locally freely (for example, via left or right translation), local freeness of the group action on $M$ is also a necessary condition for the existence of a moving frame.
\begin{Theorem}\label{local freeness and moving frames}
Let $M$ be a $G$-space and $G$ act on itself locally freely.  Then a moving frame $\rho:M\longrightarrow G$ exists if and only if the action of $G$ on $M$ is locally free.
\end{Theorem}

\begin{proof}
If $G$ acts on $M$ locally freely, then we can use the method of choosing a cross-section to the group orbits to construct a moving frame as explained above.

Conversely, suppose that $\rho:M\longrightarrow G$ is a moving frame.  Let $z\in M$ be given.  There exists an open neighborhood $V$ of the identity element $e$ in $G$ such that, for any $h\in V$, $\rho(h\cdot z)=h\cdot\rho(z)$.  Now, with respect to the action of $G$ on itself, which is a locally free action, the stabilizer of $\rho(z)\in G$, $G_{\rho(z)}:=\{ g\in G \, | \, g\cdot \rho(z)=\rho(z)  \}$, is discrete, and thus there is an open neighborhood $W$ of $e$ in $G$ such that $W\cap G_{\rho(z)}=\{e\}$.  Let $G_z:=\{ h\in G \, | \, h\cdot z=z\}$ be the stabilizer of $z$.  If $h\in (V\cap W)\cap G_z$, then $h\cdot\rho(z)=\rho(h\cdot z)=\rho(z)$, implying that $h\in W\cap G_{\rho(z)}$, and thus $h=e$.  This means that $G_z$ is discrete.  Therefore, the action of $G$ on $M$ is locally free.
\end{proof}

\begin{Remark}
We will adhere to the usual method of using cross-sections to construct moving frames, and $G$ is assumed acting on itself by the rule $g\cdot h\longmapsto hg^{-1}$ as all our examples will involve only {\em left} $G$-actions on $M$.  Thus when it comes to discussing $G$-invariant dif\/ferential forms on~$G$, such as Maurer--Cartan forms, their invariance shall be meant to be {\em right}-invariance unless explicitly specif\/ied otherwise.
\end{Remark}

Now we continue with the previous example and demonstrate how to construct a moving frame.

\begin{Example}[Euclidean plane]
Recall from Example~\ref{SE(2) action original} the prolonged action $\mathcal{A}:{\rm SE}(2)\times J^1(\R,\R^2)\longrightarrow J^1(\R,\R^2)$ of the group $G:={\rm SE}(2)$ on $M:=J^1(\R,\R^2)$ in local coordinates: $(\theta,a,b)\cdot (x,u,v,u_x,v_x)=(x,\widehat{u},\widehat{v},\widehat{u_x}, \widehat{v_x})$.
We choose the cross-section to the ${\rm SE}(2)$-orbits in $J^1(\R,\R^2)$, characterized by the  (normalization) equations:
\begin{gather}\label{SE(2) normalization}
\mathcal{A}^*u=\widehat{u}=0, \qquad \mathcal{A}^*v=\widehat{v}=0, \qquad  \mathcal{A}^*u_x=\widehat{u_x}=0.
\end{gather}
These equations are solved to f\/ind the group parameters
\[
\theta=\tan^{-1}\left(\frac{u_x}{v_x}\right), \qquad a=\frac{vu_x-uv_x}{\sqrt{u_x^2+v_x^2}}, \qquad b=\frac{-uu_x-vv_x}{\sqrt{u_x^2+v_x^2}},
\]
that uniquely determine a group element in a suf\/f\/iciently small neighborhood of the identity element of ${\rm SE}(2)$.  Then the corresponding moving frame $\rho: J^1(\R, \R^2)\longrightarrow {\rm SE}(2)$ is given by
\[
\rho(x,u,v,u_x,v_x)=\left(\tan^{-1}\left(\frac{u_x}{v_x}\right),\frac{vu_x-uv_x}{\sqrt{u_x^2+v_x^2}},\frac{-uu_x-vv_x}{\sqrt{u_x^2+v_x^2}}\right)
\]
and its moving frame section $\sigma:J^1(\R, \R^2)\longrightarrow {\rm SE}(2)\times J^1(\R, \R^2)$ is
\[
\sigma(x,u,v,u_x,v_x)=\bigl(\rho(x,u,v,u_x,v_x), x,u,v,u_x,v_x\bigr).
\]
The existence of a moving frame implies that ${\rm SE}(2)$ acts on $J^1(\R, \R^2)$ locally freely according to Theorem~\ref{local freeness and moving frames}.
\end{Example}
In general, let us denote the action of $G$ on $M$ by $\Act:G\times M\longrightarrow M$, $(g,z)\longmapsto g\cdot z:=\Act(g,z)$, and def\/ine an action of $G$ on $G\times M$ by
\[
G\times (G\times M)\longrightarrow G\times M, \qquad (g,(h,z))\longmapsto (g\cdot h,g\cdot z),
\]
where $g\cdot h$ ref\/lects any locally free action of $G$ on itself such as $g\cdot h\longmapsto hg^{-1}$.  Suppose that $\rho:M\longrightarrow G$ is a moving frame, $\sigma:M\longrightarrow G\times M$ its associated moving frame section, and $\pi:G\times M\longrightarrow M$ the canonical projection.

\begin{Proposition}\label{various maps}
The maps $\sigma$ and $\pi$ are $G$-equivariant, and $\Act$ is $G$-invariant.  Also the map $\iota:=\Act\circ \sigma:M\longrightarrow M$ is $G$-invariant.
\end{Proposition}

\begin{proof}
We will show that the moving frame section $\sigma:M\longrightarrow G\times M$ is equivariant.  The other claims about $\pi$, $\Act$, and $\iota$ can all be proved in a similar fashion.  Thus suppose that $g\in G$ and $z\in M$.  Then $\sigma(g\cdot z)=(\rho(g\cdot z),g\cdot z)=(g\cdot\rho(z),g\cdot z)=g\cdot (\rho(z),z)=g\cdot\sigma(z)$.  Therefore $\sigma$ is $G$-equivariant.
\end{proof}

All these maps are put together in the following diagram for reference.
\[
\xymatrix{& G\times M \ar[dl]^{\pi} \ar[dr]^{\Act} & \\
          M \ar@/^/[ur]^{\sigma}  \ar[rr]_{\iota} & & M  }
\]
The {\em invariantization operator} $\iota^*:\Omega^{\bullet}(M)\longrightarrow\Omega^{\bullet}(M)$, where $\Omega^\bullet(M):=\bigoplus_{k=0}^{\dim M}\Omega^k(M)$ denotes the exterior algebra of dif\/ferential forms on $M$, turns, thanks to Proposition~\ref{equivariant-invariant}, any dif\/ferential form on $M$ into an invariant dif\/ferential form, and thus occupies a central position in many works in the literature on the theory and applications of equivariant moving frames.

\begin{Example}[Euclidean plane]
The ${\rm SE}(2)$-action def\/ined in Example~\ref{SE(2) action original} on $J^1(\R,\R^2)$, 
\[
\Act: \ {\rm SE}(2)\times J^1(\R,\R^2) \longrightarrow J^1(\R,\R^2),
\]
pulls back the coordinate functions $(x,u,v,u_x,v_x)$ of $J^1(\R,\R^2)$ to ${\rm SE}(2)$-invariant functions
\begin{gather*}
\Act^*x=x, \qquad \Act^*u=u\cos\theta-v\sin\theta+a, \qquad \Act^*v=u\sin\theta+v\cos\theta+b,
\\
\Act^*u_x=u_x\cos\theta-v_x\sin\theta, \qquad \Act^*v_x=u_x\sin\theta+v_x\cos\theta
\end{gather*}
on ${\rm SE}(2)\times J^1(\R,\R^2)$.  Pulling these functions further back by the moving frame section $\sigma:J^1(\R,\R^2)\longrightarrow {\rm SE}(2)\times J^1(\R,\R^2)$ results in ${\rm SE}(2)$-invariant functions on $J^1(\R,\R^2)$:
\begin{gather*}
\sigma^*\mathcal{A}^*x=\iota^*x=x, \qquad
\sigma^*\mathcal{A}^*u=\iota^*u=0, \qquad
\sigma^*\mathcal{A}^*v=\iota^*v=0, \\
\sigma^*\mathcal{A}^*{u_x}=\iota^*{u_x}=0, \qquad
\sigma^*\mathcal{A}^*{v_x}=\iota^*{v_x}=\sqrt{u_x^2+v_x^2}.
\end{gather*}
In connection with certain coframe f\/ields that we will construct later on (Theorem~\ref{constant-structure coframe thm}), we make the observation that the invariantization of the coordinate functions $u$, $v$, $u_x$ that were used in setting up the normalization equations~(\ref{SE(2) normalization}), necessarily and trivially, yields  {\em constants} (same as the ones put on the right-hand sides of the normalization equations).
\end{Example}

\section{Congruence of immersions}

In this section, we give an answer to our key question: what are the conditions that  immersions must satisfy in order for them to be congruent under the action of a Lie group on their ambient space?

First of all, the very def\/inition of prolongation of group actions implies that one can replace a given congruence problem of immersions by a {\em prolonged congruence problem} without changing f\/inal solutions.  More specif\/ically,
\begin{Lemma}\label{prolonged congruence}
Two immersions $\psi, \phi: X\longrightarrow M$ are congruent at $x_0\in X$ if and only if, for each $k=0,1, \dots $, their prolonged graphs $j^k\gamma(\psi), j^k\gamma(\phi): X\longrightarrow \JkXM$ are locally congruent at~$x_0$ where $G$ acts on $\JkXM$ by verticalized prolongation.
\end{Lemma}

Suppose that $\{ \omega^1, \omega^2, \dots , \omega^n\}$ is a coframe f\/ield on $M$.  The coef\/f\/icients $f^i_{jk}$ in their structure equations
\[
\dd\omega^i=-\frac{1}{2}\sum_{j,k=1}^nf^i_{jk}\omega^j\wedge\omega^k, \qquad i=1,2, \dots , n,
\]
where $f^i_{jk}+f^i_{kj}=0$, $j, k=1,2,\dots, n$,
are, in general, non-constant functions on~$M$.  They do become constants if, for example, the coframe f\/ield is the Maurer--Cartan coframe f\/ield on a Lie group or its generalization on a homogeneous space.  Our solution to the congruence problem will rely on a particular kind of invariant coframe f\/ields that have {\em constant} structure functions.

\begin{Definition}\label{define constant-structure coframes}
A coframe f\/ield on $M$ is said to be of {\em constant structure} if all its structure functions are constant functions.
\end{Definition}

The following theorem, one of our main results of the paper, f\/inds a {\em necessary and sufficient} condition for immersions to be congruent when there exists on the ambient $G$-space a $G$-invariant coframe f\/ield of constant structure.
\begin{Theorem}\label{general solution}
Suppose that an $n$-dimensional $G$-space $M$ has a $G$-invariant coframe field~$\{\omega^i\}$ of constant structure.  Let $\psi, \phi:X\longrightarrow M$ be immersions.  Then $\psi$ and $\phi$ are congruent at $x_0\in X$ if and only if there exist $g\in G$ and $x_0\in X$ such that $\psi(x_0)=g\circ\phi(x_0)$ and $\psi^*\omega^i=\phi^*\omega^i$, $i=1,2,\dots ,n$, on a neighborhood of~$x_0$.
\end{Theorem}

\begin{proof}
The proof for the ``only if'' direction is trivial.  To prove the ``if" direction, let $\pi_X$ and~$\pi_M$ denote the canonical projections of $X\times M$ onto the factors $X$ and $M$, respectively.   Suppose that $\{C^i_{jk}\}$ are the structure constants of the coframe f\/ield~$\{\omega^i\}$ so that for each $i=1,2,\dots , n$,
\[
\dd\omega^i=-\frac{1}{2}\sum_{j,k=1}^nC^i_{jk}\omega^j\wedge\omega^k,
\]
 where $C^i_{jk}+C^i_{kj}=0$, $j,k=1,2,\dots , n$.
The maps that we are looking at are shown in the following diagram
\[
\xymatrix{& X\times M \ar[dl]_{\pi_X} \ar[dr]^{\pi_M} & \\
          X  \ar[rr]_{\psi} & & M  }
\]
Consider the ideal
$\mathcal{I}:=\langle \pi_M^*\omega^i-\pi_X^*\psi^*\omega^i  \rangle \subset \Omega^\bullet(X\times M)$
generated algebraically by the one-forms $\{ \pi_M^*\omega^i-\pi_X^*\psi^*\omega^i \}$.   Then, for each $i=1,2, \dots ,n$,
\begin{gather*}
 \dd \big(\pi_M^*\omega^i-\pi_X^*\psi^*\omega^i\big)
= -\frac{1}{2}\sum_{j,k=1}^nC^i_{jk}\pi_M^*\omega^j\wedge\pi_M^*\omega^k+\frac{1}{2}\sum_{j,k=1}^nC^i_{jk}\pi_X^*\psi^*\omega^j\wedge\pi_X^*\psi^*\omega^k\\
\hphantom{\dd \big(\pi_M^*\omega^i-\pi_X^*\psi^*\omega^i\big) }{}
= -\frac{1}{2}\sum_{j,k=1}^nC^i_{jk}\pi_M^*\omega^j\wedge\pi_M^*\omega^k+\frac{1}{2}\sum_{j,k=1}^nC^i_{jk}\pi_M^*\omega^j\wedge\pi_X^*\psi^*\omega^k\\
\hphantom{\dd \big(\pi_M^*\omega^i-\pi_X^*\psi^*\omega^i\big)= }{}
-\frac{1}{2}\sum_{j,k=1}^nC^i_{jk}\pi_M^*\omega^j\wedge\pi_X^*\psi^*\omega^k  +\frac{1}{2}\sum_{j,k=1}^nC^i_{jk}\pi_X^*\psi^*\omega^j\wedge\pi_X^*\psi^*\omega^k\\
\hphantom{\dd \big(\pi_M^*\omega^i-\pi_X^*\psi^*\omega^i\big) }{}
= -\frac{1}{2}\sum_{j,k=1}^nC^i_{jk}\pi_M^*\omega^j\wedge\big(\pi_M^*\omega^k-\pi_X^*\psi^*\omega^k\big)\\
\hphantom{\dd \big(\pi_M^*\omega^i-\pi_X^*\psi^*\omega^i\big)= }{}
-\frac{1}{2}\sum_{j,k=1}^nC^i_{jk}\big(\pi_M^*\omega^j-\pi_X^*\psi^*\omega^j\big)\wedge\pi_X^*\psi^*\omega^k
\equiv   0  \mod  \mathcal{I}.
\end{gather*}
This shows that $\mathcal{I}$ is a {\em differential} ideal.   Since, by assumption, $g\circ\phi(x_0)=\psi(x_0)$ for some $g\in G$ and $x_0\in X$, and since  $(g\circ\phi)^*\omega^i=\phi^*g^*\omega^i=\phi^*\omega^i=\psi^*\omega^i$, the uniqueness result in {\em Cartan's technique of the graph}, or Theorem~2.34 in~\cite{warner}, implies that $g\circ\phi(x)=\psi(x)$ for all $x$ in an open neighborhood of $x_0$.  This completes the proof.
\end{proof}

Therefore this theorem essentially turns the congruence problem of immersions into the problem of constructing an invariant coframe f\/ield of constant structure.  In general, of course, there is no guarantee that a given $G$-space ~$M$ will admit any invariant coframe f\/ield, let alone one of constant structure.  However, in many situations, a {\em prolonged} verticalized action of~$G$ on a jet bundle $\JkXM$ may become locally free for some suf\/f\/iciently large order~$k$, which will then lead to the construction of a~$\kth$ order moving frame and hence an invariant coframe f\/ield on the jet bundle (to be explained below).  Once we have in particular a {\em constant-structure} invariant coframe f\/ield constructed on $\JkXM$, we can solve, thanks to Lemma~\ref{prolonged congruence}, the original congruence problem of immersions $\psi: X\longrightarrow M$ by solving, instead, the congruence problem of the prolonged graphs $j^k\gamma(\psi): X\longrightarrow \JkXM$  of the immersions $\psi$.  The following diagram depicts how the immersion $\psi$ in the congruence problem is converted to a prolonged immersion~$j^k\gamma (\psi)$  of its graph~$\gamma(\psi)$ so that we can try our method of using a constant-structure invariant coframe f\/ield on~$\JkXM$
\[
\xymatrix{& & \\
X  \ar[rr]^{\psi} & &  M  }
\xymatrix@C=8em@R=1em{  \ar@/^/[r]  &   }
\xymatrix{& &  \JkXM \ar[d]^{}\\
          X  \ar[rru]^{j^k\gamma(\psi)} \ar[rr]_{\gamma(\psi)} & & X\times M  }
\]

\begin{Remark}\label{griffiths}
We of\/fer a brief perspective, for the reader familiar with the inceptive paper~\cite{griffiths}, on how our result compares to some part of that work.  Grif\/f\/iths in \cite{griffiths} bases his foundational lemmas determining congruence, proved in the paper's f\/irst section, upon the crucial transitivity assumption that any two smooth maps $f$ and $\tilde{f}$ into a Lie group~$G$ (or a homogeneous $G$-space) are related by a $G$-valued function~$h$ so that $f(x)=h(x)\tilde{f}(x)$ for all~$x$ in the domain of the maps.  This transitivity assumption is conspicuously absent in our take on the problem, and thus Grif\/f\/iths' result is subsumed as a special case of our more general result, Theorem~\ref{general solution}.  On the other hand, our constant-structure invariant coframe f\/ields can be interpreted as a natural extension to arbitrary $G$-spaces, homogeneous or not, of Maurer--Cartan coframes (which are always of constant structure) that Grif\/f\/iths uses to gauge the dif\/ference $h$ between the two maps~$f$ and~$\tilde{f}$.
\end{Remark}

While Theorem~\ref{general solution} does not assume transitivity of the action of $G$ on $M$, and thus can solve the congruence problem of immersions in both homogeneous and nonhomogeneous spaces within one coherent framework, its method requires that the invariant coframe f\/ield be of constant structure.  Therefore, we proceed to discuss the existence and construction of such a coframe f\/ield.  We should note that the basic ideas in the following lemmas have originated from the  references \cite{fels-olver1,fels-olver2}.

Let $\pi_G:G\times M\longrightarrow G$ denote the canonical projection $(g,z)\longmapsto g$.  Let $(u^1, \dots, u^n)$ be the local coordinate system of $M$, and $\{\mu^j \in \Omega^1(G) \, | \, j=1, \dots ,r\}$ the Maurer--Cartan coframe f\/ield of $G$.

\begin{Lemma}
The pulled-back one-forms{\samepage
\begin{gather} \label{lifted coframe}
\big\{\mathcal{A}^*\dd u^i, \pi^*_G\mu^j \in \Omega^1(G\times M) \, | \, i=1,\dots ,n; \ j=1,\dots ,r\big\}
 \end{gather}
form a $G$-invariant coframe field of $G\times M$.}
\end{Lemma}

\begin{proof}
The map $\mathcal{A}$ is $G$-invariant and $\pi_G$ is $G$-equivariant. (Recall Proposition~\ref{various maps}.)  Thus, by Proposition~\ref{equivariant-invariant}, the one-forms (\ref{lifted coframe}) are $G$-invariant.

To prove that (\ref{lifted coframe}) is a coframe f\/ield of $G\times M$, we f\/irst note that $\mathcal{A}$ is a submersion since, for any $z\in M$ and $v\in T_zM$, we can take $w:=(0,v)\in T_eG\oplus T_zM\cong T_{(e,z)}(G\times M)$ so that $\dd\mathcal{A}(w)=v$.  Thus, $\{\mathcal{A}^*\dd u^i \, | \, i=1,\dots ,n\}$ are linearly independent at every point on $G\times M$.  Likewise, $\{\pi_G^*\mu^j \, | \, j=1,\dots ,r\}$ are linearly independent on $G\times M$ since $\pi_G$ is a submersion as well.  To show that the union $\{\mathcal{A}^*\dd u^i, \pi_G^*\mu^j \, | \, i=1,\dots ,n; \ j=1,\dots ,r\}$ is also  linearly independent, it suf\/f\/ices to show that, for any $\xi\in T^*G$, the equation $\pi_G^*\xi=\sum_{i=1}^nc_i\mathcal{A}^*\dd u^i$, $c_i\in \R$, forces $c_i=0$ for all $i=1,\dots , n$.  Indeed, for any given $z\in M$, take any $w=(0,v)\in  T_{(e,z)}(G\times M)$ so that $\dd\pi_G(w)=0$.  Then the equation $(\pi_G^*\xi)(w)=\Big(\sum\limits_{i=1}^nc_i\mathcal{A}^*\dd u^i\Big)(w)$ becomes $0=\xi(\dd\pi_G(w))=\sum\limits_{i=1}^nc_i\dd u^i(\dd\mathcal{A}(w))=\sum\limits_{i=1}^nc_i\dd u^i(v)$.  Since this equation holds for all $v\in T_zM$ and all $z\in M$, it reduces to $0=\sum\limits_{i=1}^nc_i\dd u^i$ on $M$.  Thus $c_i=0$ for all $i=1,\dots, n$. Therefore  $\{\mathcal{A}^*\dd u^i, \pi_G^*\mu^j\}$ are linearly independent at every point of $G\times M$.  Finally, since $\dim (G\times M)=r+n$, the set $\{\mathcal{A}^*\dd u^i, \pi_G^*\mu^j \, | \, i=1,\dots ,n; \ j=1,\dots ,r\}$ must be a coframe f\/ield of $G\times M$.
\end{proof}

In the following is shown the prime role in our paper of moving frames and their associated sections, which is to construct constant-structure invariant coframe f\/ields on $M$

\begin{Lemma}\label{constant-structure coframe}
If the $G$-action on $M$ is locally free, then there exists on $M$ an invariant coframe field of constant structure.
\end{Lemma}
\begin{proof}
Suppose that $\sigma:M\longrightarrow G\times M$ is a moving frame section whose existence is guaranteed by Theorem~\ref{local freeness and moving frames}.  We pull back the $G$-invariant coframe f\/ield (\ref{lifted coframe}) of $G\times M$ by $\sigma$ to $M$ to obtain
\begin{gather} \label{constant-structure coframe redundant}
\sigma^*\mathcal{A}^*\dd u^i, \quad i=1,\dots, n, \qquad \sigma^*\pi_G^*\mu^j, \quad j=1,\dots ,r.
 \end{gather}
First, note that, since $\sigma$ is an embedding, the one-forms (\ref{constant-structure coframe redundant}), when evaluated at any point $z\in M$, should span the cotangent space $T^*_zM$.  Thus  $n=\dim M$ of the $n+r$ one-forms (\ref{constant-structure coframe redundant}) will form a~(local) coframe f\/ield of $M$.  Also, $G$-invariance of (\ref{constant-structure coframe redundant}) follows from the fact that the one-forms are the pull-backs of the invariant one-forms (\ref{lifted coframe}) by the $G$-equivariant map $\sigma$.  Finally, the invariant coframe f\/ield is of constant structure since $\dd(\sigma^*\mathcal{A}^*\dd u^i)=\dd\dd\sigma^*\mathcal{A}^*u^i=0$, $i=1,\dots ,n$, and $\dd(\sigma^*\pi_G^*\mu^j)=\sigma^*\pi_G^*\dd\mu^j$, $j=1,\dots , r$, where the Maurer--Cartan coframe f\/ield $\{\mu^j\}$ of $G$ is of constant structure.
\end{proof}

Recall that, for our usual initial step in constructing a moving frame, we specify a cross-section $\Gamma$ to the $G$-orbits in $M$ by setting up $r=\dim G$ normalization equations $\mathcal{A}^*u^{i_\kappa}=c_\kappa$, $\kappa=1,\dots , r$, where we typically choose to use {\em constants} for~$c_\kappa$.  An implication of the constants~$c_\kappa$ is that~$r$ of the one-forms~(\ref{constant-structure coframe redundant}) will necessarily vanish for $\sigma^*\mathcal{A}^*\dd u^{i_\kappa}=\dd c_\kappa=0$.  Therefore this typical approach gives rise to the following f\/inal formula for constructing a constant-structure invariant coframe f\/ield of~$M$.  We continue to assume that~$G$ acts on~$M$ locally freely and that~$u^i$, $i=1,\dots ,n$, are the local coordinate functions of~$M$.

\begin{Theorem}\label{constant-structure coframe thm}
Suppose that, in the course of constructing a moving frame, the $r$ coordinate functions $u^i$, $i=1,\dots ,r$, of $M$ are normalized to constants; that is, $\iota^*u^i=\sigma^*\mathcal{A}^*u^i=c_i$,  constant.  Then $M$ has the following constant-structure $G$-invariant coframe field:
\begin{gather} \label{constant-structure coframe formula}
\sigma^*\mathcal{A}^*\dd u^i=\dd\sigma^*\mathcal{A}^*u^i=\dd\iota^*u^i , \quad i=r+1,\dots, n,  \qquad \sigma^*\pi_G^*\mu^j=\rho^*\mu^j, \quad j=1,\dots ,r.
\end{gather}
\end{Theorem}

\begin{Remark}\label{not really redundant invariant coframe field}
Note that, if a coordinate function $u^\kappa$ of $M$ is normalized to a constant $c_\kappa$ in the process of choosing a cross-section $\Gamma$ to the group orbits, then $\dd\iota^*u^\kappa=\dd c_\kappa=0$.  Thus, even if we implement the formula (\ref{constant-structure coframe redundant}) for {\em every} coordinate function $u^i$, $i=1,\dots ,n$, it will not cause any problem in constructing our desired constant-structure invariant coframe f\/ield, apart from the potential unexciting prospect of having to spend some extra time computing for those vanishing phantom dif\/ferential forms.
\end{Remark}

Now suppose that  $(x^1,x^2,\dots ,x^p)$ are local coordinate functions on $X$ and that a constant-structure invariant coframe f\/ield $\Big\{\omega^i =\sum\limits_{j=1}^n f^i_j \dd u^j \, | \, i=1,2,\dots , n\Big\}$ exists on $M$ with some functions $f^i_j:M\longrightarrow \R$.  Then, assuming $p=\dim X$, the pull-backs
\begin{gather*} 
\psi^*\omega^i=\sum_{l=1}^{p}\sum\limits_{j=1}^n(f^i_j\circ\psi)\frac{\partial(u^j\circ\psi)}{\partial x^l}\dd x^l, \qquad i=1,2,\dots , n,
 \end{gather*}
by congruent immersions $\psi:X\longrightarrow M$ must all agree by Theorem~\ref{general solution}.
In turn, these pulled-back forms will agree if and only if the dif\/ferential functions
$\Big\{\sum\limits_{j=1}^n f^i_j u^j_l\Big\}$, where $u^j_l$ are the local coordinate functions on the f\/ibers of the bundle $J^1(X, M)\longrightarrow X\times M$, agree when pulled back by the immersions $\psi$.  Thus we have proved:
\begin{Lemma}\label{complete set of invariants}
The differential functions $\Big\{\sum\limits_{j=1}^nf^i_j u^j_l \, | \, i=1,\dots , n; \  l=1,\dots , p\Big\}$ form a {\em complete} system of congruence invariants for immersions from~$X$ to~$M$.
\end{Lemma}

\begin{Remark}
In general, the complete system $\Big\{\sum\limits_{j=1}^nf^i_j u^j_l\Big\}$ will not be functionally independent, and thus there will be some functional redundancy in the system, as is indeed the case with most of the examples in our paper.  Also, a subset consisting of only functionally independent invariants may still be further reduced to an even smaller subset consisting of lower order dif\/ferential invariants by employing the idea of {\em invariant differential operators}, \cite{equivalence-book}, that can be used to eliminate dif\/ferential redundancy.  Related to this point, it should be noted that the case of {\em variable-parameter} equivalence problems has in recent years witnessed major advances made by Olver et al.\ both in the general approach to obtaining (minimal) generating systems of dif\/ferential invariants and in its applications to some specif\/ic geometric settings; see~\cite{hubert,surfaces,centro-affine,dia}.  In our present paper, we will not venture into this intricate question of {\em minimalizing} systems of invariants; the interested reader should instead consult the cited references.
\end{Remark}

In general, a group action on $M$ may not be locally free, and so a constant-structure invariant coframe f\/ield may not be available on $M$.  However, if we prolong the verticalized action on $X\times M$ of the group to a jet bundle $\JkXM$ of a suf\/f\/iciently high order $k$, the resulting action may become locally free essentially thanks to, intuitively speaking, the {\em wider} space that mitigates the crowding of the group orbits.  Once we have a locally free action on a jet bundle $\JkXM$, we can apply the ideas of Theorems~\ref{general solution} and~\ref{constant-structure coframe thm} to construct a constant-structure invariant coframe f\/ield on $\JkXM$ and solve the congruence problem.

The following result tells us basically how many invariants at most, in terms of their dif\/fe\-ren\-tial order, we should require to resolve a congruence problem.
\begin{Corollary}\label{minimal differential order}
Suppose that a $k^{th}$ order jet bundle $\JkXM$ admits a constant-structure $G$-invariant coframe field.  Then the differential order of a minimal complete set of congruence invariants is bounded by $k+1$.
\end{Corollary}

\begin{proof}
Suppose that $\Big\{ \sum\limits_{l=1}^p h_l^\alpha \dd x^l, \ \sum\limits_{|J|=0}^k \sum\limits_{j=1}^n f^{I,J}_{i,j} \dd u^j_J \Big\}$ (where  $\alpha=1,\dots,p$,  $i=1,\dots, n$, and $I$~is a symmetric multi-index over $\{ 1, \dots ,p \}$ of order up to~$k$) is a constant-structure invariant coframe f\/ield on $\JkXM$.  Note that all the dif\/ferential functions $h^\alpha_l$, $x^l$, $f^{I,J}_{i,j}$, and $u^j_J$ are of order up to $k$.   Then, by Lemma~\ref{complete set of invariants} and with the understanding that $M$  and $\psi$ in the proof of the lemma are now replaced respectively by $\JkXM$ and $j^k\gamma(\psi)$, a minimal complete set of congruence invariants must be a subset of $\Big\{ h^\alpha_l, \sum\limits_{|J|=0}^k \sum\limits_{j=1}^n f^{I,J}_{i,j} u^j_{J,l} \Big\}$ which is of dif\/ferential order $k+1$ or less.
\end{proof}

Thus, for example, if the verticalized prolongation of the group action becomes locally free at order $k$, then the conclusion of the corollary holds true.

\begin{Remark}
Note that f\/ixed-parameter congruence is a stronger condition than variable-parameter cong\-ruen\-ce.  Thus, Corollary~\ref{minimal differential order} implies in the case of variable-parameter congruen\-ce problems that, under the assumptions of the corollary, the dif\/ferential order of a minimal complete set of {\em variable-parameter} congruence invariants is also bounded by $k+1$.
\end{Remark}

Once a constant-structure invariant coframe f\/ield $\{\omega^i\}$ exists, any other invariant coframe f\/ield $\{\zeta^i\}$ of the form
\[
\zeta^i:=\sum_{j=1}^n a^i_j\omega^j, \qquad a^i_j= \text{const}, \qquad i=1,2, \dots , n,
\]
will be of constant structure, and can also be used for solving the congruence problem of immersions.  If, however, an invariant coframe f\/ield $\{\zeta^i\}$ is not of constant structure, then at least {\em necessary} conditions for congruence can be obtained.

\begin{Proposition}\label{non-constant structure prop}
Suppose that $M$ admits an invariant coframe field of constant structure and that $\{\zeta^i\}$ is any invariant coframe field that may not be of constant structure.  If two immersions $\psi, \phi:X\longrightarrow M$ are congruent, then $\psi^*\zeta^i=\phi^*\zeta^i$ for all $i=1,2,\dots , n$.
\end{Proposition}

\begin{proof}
Let $\{\omega^i\}$ be a constant-structure invariant coframe f\/ield on $M$.  Suppose $\psi=g\circ\phi$ for some $g\in G$.
Since $\{\zeta^i\}$ is an invariant coframe f\/ield, there exist invariant functions $f^i_j$ on $M$ such that $\zeta^i=\sum\limits_{j=1}^n f^i_j\omega^j$ for all $i=1, \dots , n$.  Then for each $i$,
\begin{gather*}
\psi^*\zeta^i=\sum_{j=1}^n \!\big((g\circ\phi)^*f^i_j\big)\big((g\circ\phi)^*\omega^j\big)=\sum_{j=1}^n\!\big(\phi^*g^*f^i_j\big)\big(\phi^*g^*\omega^j\big)
=\sum_{j=1}^n\!\big(\phi^*f^i_j\big)\big(\phi^*\omega^j\big)=\phi^*\zeta^i.\!\!\!\!\!\!\!\tag*{\qed}
\end{gather*}
\renewcommand{\qed}{}
\end{proof}

\begin{Remark}
Fortunately, we have the general formula (\ref{constant-structure coframe formula}) for constructing {\em constant-structure} invariant coframe f\/ields that will enable us to solve congruence problems completely and obtain {\em necessary and sufficient} conditions for congruence.
\end{Remark}

\begin{Example}[Euclidean plane]
We continue the discussion and notation of the ${\rm SE}(2)$-action that we had in the earlier examples.  To f\/ind a complete system of ${\rm SE}(2)$-congruence invariants of curve immersions in Euclidean plane, we need to  construct a constant-structure ${\rm SE}(2)$-invariant coframe f\/ield on $J^1(\R,\R^2)$ where we conf\/irmed that ${\rm SE}(2)$ acts locally freely.  We follow the formula (\ref{constant-structure coframe formula}) to construct such a coframe f\/ield.  According to the formula, part of the coframe f\/ield comes from the Maurer--Cartan coframe f\/ield of ${\rm SE}(2)$, which we f\/ind by f\/irst embedding ${\rm SE}(2)$ in ${\rm GL}(3)$:
\[
{\rm SE}(2)\longrightarrow {\rm GL}(3,\R), \qquad (\theta,a,b)\longmapsto g:=
\begin{bmatrix}
\cos\theta & -\sin\theta & a \\
\sin\theta & \cos\theta & b \\
0 & 0 & 1
\end{bmatrix},
\]
which acts on the hyperplane $\{ [u \ v \ 1]^{\rm T}\in\R^3 \, | \, u,v \in\R \}$ via matrix multiplication,
and then taking the linearly independent entries of the $\mathfrak{gl}(3,\R)$-valued Maurer--Cartan form $(\dd g) g^{-1}$:
\[ 
\mu^1:=\dd\theta, \qquad \mu^2:=\dd a+b\dd\theta, \qquad \mu^3:=\dd b-a\dd\theta.
\]

We pull these Maurer--Cartan forms by the moving frame $\rho:J^1(\R,\R^2)\longrightarrow {\rm SE}(2)$ to obtain the invariant one-forms on $J^1(\R,\R^2)$:
\begin{gather}\label{SE(2) Maurer--Cartan coframe}
\rho^*\mu^1=\frac{v_x\dd u_x-u_x\dd v_x}{u_x^2+v_x^2}, \qquad \rho^*\mu^2=\frac{-v_x\dd u+u_x\dd v}{\sqrt{u_x^2+v_x^2}}, \qquad \rho^*\mu^3=\frac{-u_x\dd u-v_x\dd v}{\sqrt{u_x^2+v_x^2}},
\end{gather}
which will be part of the coframe f\/ield that we are trying to construct.

The other part of the formula (\ref{constant-structure coframe formula}) is implemented to obtain more invariant one-forms on $J^1(\R,\R^2)$:
\[ 
\dd\iota^*x=\dd x, \qquad \dd\iota^*v_x=\frac{u_x}{\sqrt{u_x^2+v_x^2}}\dd u_x+\frac{v_x}{\sqrt{u_x^2+v_x^2}}\dd v_x,
\]
that complement (\ref{SE(2) Maurer--Cartan coframe}) to form a constant-structure invariant coframe f\/ield of $J^1(\R,\R^2)$.  The f\/inal step, suggested by Theorem~\ref{general solution}, is to pull these coframe forms by the prolonged graph $j^1\gamma(\psi)$ of a generic curve immersion  $\psi:\R\longrightarrow \R^2$, $x\longmapsto (u(x),v(x))$:
\begin{gather}
\big(j^1\gamma(\psi)\big)^*\mu^1=\frac{u_{xx}v_x-u_xv_{xx}}{u_x^2+v_x^2}\dd x, \qquad
\big(j^1\gamma(\psi)\big)^*\mu^2=0,\nonumber \\
\big(j^1\gamma(\psi)\big)^*\mu^3=-\sqrt{u_x^2+v_x^2}\dd x, \qquad
\big(j^1\gamma(\psi)\big)^*\dd\iota^*x=\dd x, \label{SE(2) final pullback}\\
\big(j^1\gamma(\psi)\big)^*\dd\iota^*v_x=\frac{u_xu_{xx}+v_xv_{xx}}{\sqrt{u_x^2+v_x^2}}\dd x.\nonumber
\end{gather}
Therefore, the coef\/f\/icient functions of $\dd x$ in (\ref{SE(2) final pullback}),
\begin{gather} \label{SE(2) complete invariants}
\frac{u_{xx}v_x-u_xv_{xx}}{u_x^2+v_x^2}, \qquad  -\sqrt{u_x^2+v_x^2}, \qquad \frac{u_xu_{xx}+v_xv_{xx}}{\sqrt{u_x^2+v_x^2}},
 \end{gather}
constitute a complete system of ${\rm SE}(2)$-invariants for curve immersions in $\R^2$.  According to Theorem~\ref{general solution}, two immersions $\R\longrightarrow \R^2$ are congruent if and only if the functions (\ref{SE(2) complete invariants}), when restricted to the two (prolonged) immersed submanifolds, yield the same values at every point in an open subset of the domain $\R$ of the immersions and also there is a transformation $g\in {\rm SE}(2)$ that sends some f\/irst-order jet $(\widetilde{x},\widetilde{u},\widetilde{v},\widetilde{u_x},\widetilde{v_x})\in J^1(\R,\R^2)$ of one of the curves to a jet of the other.

Note that the f\/irst two congruence invariants   $ \frac{u_{xx}v_x-u_xv_{xx}}{u_x^2+v_x^2}$  and  $ -\sqrt{u_x^2+v_x^2}$   generate algebraically a space of functions equivalent to the space generated by the well-known {\em curvature} and {\em speed} invariants of Euclidean planar curves.  The third invariant $\frac{u_xu_{xx}+v_xv_{xx}}{\sqrt{u_x^2+v_x^2}}$  is a dif\/ferential consequence of the second invariant, and hence does not provide any new information.
\end{Example}

\section{More examples}
\subsection[A homogeneous space:  $\mathbb{CP}^1$ under the action of ${\rm PSL}(2,\mathbb{C})$]{A homogeneous space:  $\boldsymbol{\mathbb{CP}^1}$ under the action of $\boldsymbol{{\rm PSL}(2,\mathbb{C})}$}

The example of planar curves used in the previous sections was for the congruence problem of immersions in the homogeneous space $(\R^2,{\rm SE}(2))$.  In this section, we deal with yet another example of a homogeneous space.

Consider the action by linear fractional (or M\"obius) transformations
\begin{gather} \label{SL2 action}
{\rm SL}(2,\mathbb{C})\times \mathbb{C}\longrightarrow \mathbb{C}, \qquad
\begin{bmatrix}
a & b \\
c & d
\end{bmatrix}
\cdot z \longmapsto \frac{az+b}{cz+d},
 \end{gather}
which is a local description of the transitive action
\[
{\rm PSL}(2,\mathbb{C})\times \mathbb{CP}^1\longrightarrow \mathbb{CP}^1, \qquad
\overline{
\begin{bmatrix}
a & b \\
c & d
\end{bmatrix}
}
\cdot
\overline{
\begin{bmatrix}
z \\ w
\end{bmatrix}
}
 \longmapsto
\overline{
\begin{bmatrix}
az+bw \\
cz+dw
\end{bmatrix}
}
,
\]
where the overlines represent the standard projections to the respective quotient spaces ${\rm PSL}(2,\!\mathbb{C})$ $={\rm SL}(2,\mathbb{C})/({\rm SL}(2,\mathbb{C})\cap \mathbb{C}^*I)$ and $\mathbb{CP}^1=(\mathbb{C}^2-\{(0,0)\})/\sim$.

We will study local congruence of holomorphic maps $\psi:\mathbb{CP}^1\longrightarrow \mathbb{CP}^1$ with nonvanishing derivatives.  To do so, we f\/irst rewrite the action (\ref{SL2 action}) using the parametrization of ${\rm SL}(2,\mathbb{C})$ by $(a,b,c)$:
\[ 
{\rm SL}(2,\mathbb{C})\times \mathbb{C}\longrightarrow \mathbb{C}, \qquad
\begin{bmatrix}
a & b \\
c & (bc+1)/a
\end{bmatrix}
\cdot z \longmapsto \frac{a^2z+ab}{acz+bc+1},
\]
and also view the immersions in local coordinate charts of $\mathbb{CP}^1$ so that they are of the form
$\psi:\mathbb{C}\longrightarrow \mathbb{C}$.  Next, we take the graph of $\psi$,
\[
\gamma(\psi): \ \mathbb{C}\longrightarrow \mathbb{C}\times \mathbb{C}, \qquad z\longmapsto (z,\psi(z)),
\]
where  ${\rm SL}(2,\mathbb{C})$ acts on $\mathbb{C}\times \mathbb{C}$ by verticalization:
\[ 
{\rm SL}(2,\mathbb{C})\times (\mathbb{C}\times\mathbb{C})\longrightarrow\mathbb{C}\times\mathbb{C}, \qquad
(a,b,c) \cdot (w,z) \longmapsto \left(w,\frac{a^2z+ab}{acz+bc+1}\right).
\]
This action is not locally free, which is obvious from consideration of the dimensions of ${\rm PSL}(2,\mathbb{C})$ and $\mathbb{C}\times\mathbb{C}$, and thus we need to try prolonging the action to jet bundles over $\mathbb{C}\times\mathbb{C}$ and see if the resulting action is locally free there.  In fact, the prolonged action of ${\rm SL}(2,\mathbb{C})$ on $J^2(\mathbb{C},\mathbb{C})$ turns out to be locally free.  Specif\/ically, the second-order prolonged action
$\mathcal{A}:{\rm SL}(2,\mathbb{C})\times J^2(\mathbb{C},\mathbb{C})\longrightarrow J^2(\mathbb{C},\mathbb{C})$
is given by
\begin{gather*}
 (a,b,c)\cdot (w,z,z_w,z_{ww}) \\
  \qquad \longmapsto \left(w,\frac{a^2z+ab}{acz+bc+1}, \frac{a^2 z_w}{(acz+bc+1)^2},\frac{-2a^3cz_w^2+(a^2+a^2bc)z_{ww}+a^3c z z_{ww}}{(acz+bc+1)^3}\right)
\end{gather*}
where $(w,z,z_w,z_{ww})$ denotes the local coordinate system on $J^2(\mathbb{C},\mathbb{C})$.
At this point, the chosen (normalization) equations
\begin{gather}
\mathcal{A}^*z=\frac{a^2z+ab}{acz+bc+1}=0, \qquad \mathcal{A}^*z_w=\frac{a^2 z_w}{(acz+bc+1)^2}=1,\nonumber\\
 \mathcal{A}^*z_{ww}=\frac{-2a^3cz_w^2+(a^2+a^2bc)z_{ww}+a^3c z z_{ww}}{(acz+bc+1)^3}=0,\label{SL2 normalization}
\end{gather}
can be solved, where $z_w\neq 0$, for the group parameters
\begin{gather*}
a=\frac{1}{\sqrt{z_w}}, \qquad b = -\frac{z}{\sqrt{z_w}}, \qquad c= \frac{z_{ww}}{2 z_w \sqrt{z_w}},
 \end{gather*}
which implies that the group action is locally free, and def\/ines the moving frame %
\[
\rho: \ J^2(\mathbb{C},\mathbb{C})\longrightarrow {\rm SL}(2,\mathbb{C}), \qquad (w,z,z_w,z_{ww})\longmapsto (a,b,c),
\]
and the associated section
\[
\sigma: \ J^2(\mathbb{C},\mathbb{C})\longrightarrow {\rm SL}(2,\mathbb{C})\times J^2(\mathbb{C},\mathbb{C}), \qquad (w,z,z_w,z_{ww})\longmapsto (a,b,c,w,z,z_w,z_{ww}).
\]

Now we recall the formula (\ref{constant-structure coframe formula}) to construct a constant-structure ${\rm SL}(2,\mathbb{C})$-invariant coframe f\/ield on $J^2(\mathbb{C},\mathbb{C})$.  Among the coordinate functions $(w,z,z_w,z_{ww})$ of $J^2(\mathbb{C},\mathbb{C})$, $w$ was the one that was not used in the normalization equations (\ref{SL2 normalization}) and thus it leads to the invariant one-form
\begin{gather} \label{SL2 normalized coframe}
\dd\iota^*w=\dd w
 \end{gather}
on $J^2(\mathbb{C},\mathbb{C})$ as part of the coframe f\/ield that we are constructing.
The other part of the coframe f\/ield requires, according to (\ref{constant-structure coframe formula}), that we f\/irst f\/ind Maurer--Cartan forms of ${\rm SL}(2,\mathbb{C})$.  Thus, assuming
\[
g:=
\begin{bmatrix}
a & b \\
c & (bc+1)/a
\end{bmatrix} \in {\rm SL}(2,\mathbb{C}),
\]
we take the linearly independent entries of the $\mathfrak{sl}(2,\mathbb{C})$-valued Maurer--Cartan form $(\dd g) g^{-1}$ to obtain the Maurer--Cartan coframe f\/ield of ${\rm SL}(2,\mathbb{C})$:
\[ 
\mu^1:=\frac{bc+1}{a}\dd a-c\dd b, \qquad \mu^2:=-b\dd a+a\dd b  , \qquad \mu^3:=\frac{bc^2+c}{a^2}\dd a-\frac{c^2}{a}\dd b+\frac{1}{a}\dd c   .
\]
Pulling back these Maurer--Cartan forms by the moving frame $\rho:J^2(\mathbb{C},\mathbb{C})\longrightarrow {\rm SL}(2,\mathbb{C})$ yields the following invariant one-forms on $J^2(\mathbb{C},\mathbb{C})$:
\[ 
\rho^*\mu^1=\frac{z_{ww}}{2z_w^2}\dd z-\frac{1}{2z_w}\dd z_w, \qquad \rho^*\mu^2=-\frac{1}{z_w}\dd z  , \qquad \rho^*\mu^3=\frac{z_{ww}^2}{4z_w^3}\dd z-\frac{z_{ww}}{z_w^2}\dd z_w+\frac{1}{2z_w}\dd z_{ww}
\]
that, together with (\ref{SL2 normalized coframe}), form a constant-structure ${\rm SL}(2,\mathbb{C})$-invariant coframe f\/ield on $J^2(\mathbb{C},\mathbb{C})$.

Now the f\/inal step is to pull back the one-forms in the coframe by the prolonged graph $j^2\gamma(\psi)$ of a generic holomorphic immersion $\psi:\mathbb{C}\longrightarrow\mathbb{C}$, $w\longmapsto z(w)$, to obtain:
\begin{gather*} 
 \big(j^2\gamma(\psi)\big)^*\dd\iota^*w=\dd w , \qquad \big(j^2\gamma(\psi)\big)^*\rho^*\mu^1=0, \qquad \big(j^2\gamma(\psi)\big)^*\rho^*\mu^2=-\dd w , \\
 \big(j^2\gamma(\psi)\big)^*\rho^*\mu^3=\left( \frac{z_{www}}{2z_w}-\frac{3z_{ww}^2}{4z_w^2} \right)\dd w .
\end{gather*}

Thus $  \frac{z_{www}}{2z_w}-\frac{3z_{ww}^2}{4z_w^2}$ is the only non-constant invariant in this example, and  therefore, according to Theorem~\ref{general solution}, two immersions $\mathbb{C}\longrightarrow \mathbb{C}$ are congruent under linear fractional transformations~(\ref{SL2 action}) or two immersions $\mathbb{CP}^1\longrightarrow \mathbb{CP}^1$ are congruent under the action of ${\rm PSL}(2,\mathbb{C})$ if and only if the restriction of the function $\frac{z_{www}}{2z_w}-\frac{3z_{ww}^2}{4z_w^2}$ on the two (prolonged) immersed submanifolds agree and also there exists an element $g$ of the group that takes some jet $(\widetilde{w},\widetilde{z},\widetilde{z_w},\widetilde{z_{ww}})\in J^2(\mathbb{C},\mathbb{C})$ of one of the immersions to a jet of the other immersion.  The function $\frac{z_{www}}{2z_w}-\frac{3z_{ww}^2}{4z_w^2}$, therefore, is an ${\rm SL}(2,\mathbb{C})$- or ${\rm PSL}(2,\mathbb{C})$-congruence invariant, and, when multiplied by $2$, is in fact known as the {\em Schwarzian derivative}.

\subsection[A nonhomogeneous space: $\R^3$ under the action of ${\rm SO}(3)$]{A nonhomogeneous space: $\boldsymbol{\R^3}$ under the action of $\boldsymbol{{\rm SO}(3)}$}

Consider the action of ${\rm SO}(3)$ on $\R^3$ via the standard rotations.  Note that the resulting space $(\R^3,{\rm SO}(3))$ constitutes a {\em nonhomogeneous} space.  We will treat the congruence problem of orien\-table regular hypersurfaces in this space using the same technique that we used for congruence problems for homogeneous spaces.

For this example, however, we will take a somewhat implicit way of describing the group elements, as opposed to the explicit parametric descriptions of symmetry groups that we used for the other examples.  Although it is possible to describe ${\rm SO}(3)$ by parameters, such as {\em Euler angles}, still arriving at the same f\/inal solution to our congruence problem, an implicit description of ${\rm SO}(3)$ turns out not only to reduce greatly the amount of computation that would otherwise be required by explicit presence of parameters, but also to give us a better geometric understanding, allowing for deft manipulation, of a moving frame that we will construct.

Thus, any matrix element of the group ${\rm SO}(3)$ will from now on be denoted simply by the letter $R$ without referring to any parameters.  Also, in order not to clutter expressions with too many dif\/ferent symbols, we will continue the tradition of using same symbols for multiple purposes, supported at times by the customary canonical identif\/ication of a vector space such as $\R^3$ with its tangent spaces at various points, insofar as there is no danger of serious confusion;  in particular, we will use the column matrix notation $u:=[u^1 \ u^2 \ u^3]^{\rm T}$ to denote, depending on the context, the coordinate system of $\R^3$, a point as a vector in $\R^3$, or the image of an immersion $\R^2\longrightarrow \R^3$.

As always, we view any regular surface $\psi:\R^2\longrightarrow \R^3$ in terms of its graph
\[
\gamma(\psi): \ \R^2\longrightarrow \R^2\times\R^3, \qquad \big(x^1,x^2\big) \longrightarrow \big(x^1,x^2,\psi\big(x^1,x^2\big)\big),
\]
where $x:=(x^1,x^2)$ denotes the coordinate system of the domain $\R^2$.  We verticalize the action of ${\rm SO}(3)$ to the corresponding action on $\R^2\times\R^3$:
\[
{\rm SO}(3)\times \big(\R^2\times\R^3\big)\longrightarrow \R^2\times\R^3, \qquad R\cdot (x,u) \longmapsto (x, R u),
\]
where $R u$ signif\/ies matrix multiplication.  Since this action is not locally free, we try prolonging the action in an attempt to obtain a locally free action:
\[
\mathcal{A}: \ {\rm SO}(3)\times J^1\big(\R^2,\R^3\big) \longrightarrow J^1\big(\R^2,\R^3\big), \qquad R\cdot (x,u,u_1,u_2)\longmapsto (x,Ru,Ru_1,Ru_2),
\]
where $u_1:=[u^1_1 \ u^2_1 \ u^3_1]^{\rm T}$ and $u_2:=[u^1_2 \ u^2_2 \ u^3_2]^{\rm T}$ denote the f\/iber coordinates of the bundle $J^1(\R^2,\R^3)\longrightarrow \R^2\times\R^3$.
Indeed, this f\/irst prolongation of the action is locally free since we can solve, for the group element $R$, the following chosen  (normalization) equations
\begin{gather}\label{SO3 normalization eqns}
\mathcal{A}^*u_1=R\big[u^1_1 \ u^2_1 \ u^3_1\big]^{\rm T}=[\star \ \ 0  \  \ 0]^{\rm T}, \qquad
\mathcal{A}^*u_2=R\big[u^1_2 \ u^2_2 \ u^3_2\big]^{\rm T}=[\star \ \star  \ 0]^{\rm T},
\end{gather}
where $\star$ means that the entry is not normalized and hence left free.  To help describe the solution $R$ satisfying the equations (\ref{SO3 normalization eqns}), we need f\/irst def\/ine a few $\R^3$-valued vector f\/ields along surfaces:
\[ 
\nn:=\frac{u_1\times u_2}{\lVert u_1\times u_2 \rVert}, \qquad \ttt:=\frac{u_1}{\lVert u_1 \rVert}, \qquad \ttv:=\nn\times\ttt.
\]
Note that these ordered vector f\/ields form an oriented orthonormal frame f\/ield of $\R^3$, restricted to surfaces, possessing the same orientation as the standard one for $\R^3$.  If $R$ satisf\/ies the equations~(\ref{SO3 normalization eqns}), then
\begin{gather*} 
  R\nn=[0 \ 0 \ 1]^{\rm T}, \qquad R\ttt=[1 \ 0 \ 0]^{\rm T}, \qquad \text{and} \\
  R\ttv=R(\nn\times\ttt)=(R\nn)\times (R\ttt)=[0 \ 0 \ 1]^{\rm T}\times [1 \ 0 \ 0]^{\rm T}=[0 \ 1 \ 0]^{\rm T}.
\end{gather*}
Thus,
\begin{gather*}
 R^{-1}[1 \ 0 \ 0]^{\rm T}=\ttt, \qquad R^{-1}[0 \ 1 \ 0]^{\rm T}=\ttv, \qquad R^{-1}[0 \ 0 \ 1]^{\rm T}=\nn,
\end{gather*}
and therefore
\begin{gather*}
R=[\ttt \ \ttv \ \nn]^{-1}=[\ttt \ \ttv \ \nn]^{\rm T}.
\end{gather*}

This solution $R$ is what is used in our construction of a moving frame:
\begin{gather*}
\rho: \ J^1(\R^2,\R^3)\longrightarrow {\rm SO}(3), \qquad (x,u,u_1,u_2)\longmapsto [\ttt \ \ttv \ \nn]^{\rm T},
\end{gather*}
and its associated section
\begin{gather*}
\sigma: \ J^1\big(\R^2,\R^3\big)\longrightarrow {\rm SO}(3)\times J^1\big(\R^2,\R^3\big), \quad    (x,u,u_1,u_2)\longmapsto ([\ttt \ \ttv \ \nn]^{\rm T},x,u,u_1,u_2).
\end{gather*}

To construct a constant-structure ${\rm SO}(3)$-invariant coframe f\/ield on $J^1(\R^2,\R^3)$, we now  implement the formula (\ref{constant-structure coframe formula}) while having in mind Remark~\ref{not really redundant invariant coframe field} and the def\/inition $\iota:=\mathcal{A}\circ\sigma$,
\begin{gather}
 \dd\iota^*x=\dd x, \nonumber\\
 \dd\iota^*u=\dd\sigma^*(Ru)=\dd((\rho^*R)u)=\dd([\ttt \ \ttv \ \nn]^{\rm T}  u)=[\dd \ttt \ \dd \ttv \ \dd \nn]^{\rm T}  u+[\ttt \ \ttv \ \nn]^{\rm T} \dd u,\nonumber \\
 \dd\iota^*u_1=[\dd \ttt \ \dd \ttv \ \dd \nn]^{\rm T}  u_1+[\ttt \ \ttv \ \nn]^{\rm T} \dd u_1, \nonumber\\
 \dd\iota^*u_2=[\dd \ttt \ \dd \ttv \ \dd \nn]^{\rm T}  u_2+[\ttt \ \ttv \ \nn]^{\rm T} \dd u_2.  \label{SO(3) normalized coframe}
\end{gather}
We also pull back, by the moving frame $\rho$, the $\mathfrak{so}(3)$-valued Maurer--Cartan form $(\dd R) R^{-1}$ of ${\rm SO}(3)$ to $J^1(\R^2,\R^3)$:
\begin{gather}\label{SO(3) Maurer--Cartan form}
\rho^*\big((\dd R) R^{-1}\big) =\big(\dd [\ttt \ \ttv \ \nn]^{\rm T}\big) [\ttt \ \ttv \ \nn]=
\begin{bmatrix}
0 & \ttv\cdot\dd\ttt & \nn\cdot\dd\ttt \\
-\ttv\cdot\dd\ttt & 0 & \nn\cdot\dd\ttv \\
-\nn\cdot\dd\ttt & -\nn\cdot\dd\ttv & 0
\end{bmatrix}.
\end{gather}

According to Theorem~\ref{constant-structure coframe thm} (and Remark~\ref{not really redundant invariant coframe field}), the one-forms (\ref{SO(3) normalized coframe}) and (\ref{SO(3) Maurer--Cartan form}) constitute a~cons\-tant-structure ${\rm SO}(3)$-invariant coframe f\/ield on $J^1(\R^2,\R^3)$.  Therefore, by Theorem~\ref{general solution}, two surface immersions in $\R^3$ will be congruent under the action of ${\rm SO}(3)$ if and only if there exists a transformation $R\in {\rm SO}(3)$ taking some f\/irst-order jet  $(\widetilde{x},\widetilde{u},\widetilde{u_1},\widetilde{u_2})\in J^1(\R^2,\R^3)$ of one of the surfaces to a jet  $(\widetilde{x},R\widetilde{u},R\widetilde{u_1},R\widetilde{u_2})$ of the other surface and also evaluation of the coframe f\/ield~(\ref{SO(3) normalized coframe}) and~(\ref{SO(3) Maurer--Cartan form}) results in the same one-form system for both surfaces.

To shed a bit of light on the relationship of this result with the classic case of the full orientation-preserving rigid motion group ${\rm SE}(3)={\rm SO}(3)\ltimes\R^3$, suppose that all the congruence conditions are met and thus the prolonged graphs of the two surfaces are congruent over an open subset of their domain $\R^2$.  If the vector f\/ields $u$, $u_1$, $u_2$ for each of the surfaces are expressed in terms of the orthonormal frame f\/ield $\{\ttt,  \ttv, \nn\}$ of the corresponding surface, the common restricted Maurer--Cartan forms  (\ref{SO(3) Maurer--Cartan form}) make the terms $[\dd \ttt \ \dd\ttv \ \dd\nn]^{\rm T} u$, $[\dd \ttt \ \dd\ttv \ \dd\nn]^{\rm T} u_1$, $[\dd \ttt \ \dd\ttv \ \dd\nn]^{\rm T} u_2$ in (\ref{SO(3) normalized coframe}) equal for both surfaces, and thus, in particular, the surfaces must agree on the terms
\begin{gather} \label{equivalents to fundamental forms}
[\ttt \ \ttv \ \nn]^{\rm T}\dd u, \qquad [\ttt \ \ttv \ \nn]^{\rm T}\dd u_1,  \qquad [\ttt \ \ttv \ \nn]^{\rm T}\dd u_2.
 \end{gather}
Notice that specifying these one-forms (\ref{equivalents to fundamental forms}) is equivalent (modulo some redundancy) to determining {\em the first and second fundamental forms}
\begin{gather*}
{\rm I} =(\ttt\cdot\dd u)\otimes (\ttt\cdot\dd u)+(\ttv\cdot\dd u)\otimes(\ttv\cdot\dd u), \qquad {\rm II} =\sum_{i,j=1}^2(\nn\cdot u_{ij})\dd x^i\otimes\dd x^j
 \end{gather*}
that are classically known to form a complete system of congruence invariants for surface immersions under the transitive action of ${\rm SE}(3)$.  Obviously, in view of the explanation given in the previous paragraph, the fundamental form conditions alone will not be suf\/f\/icient for congruence under the intransitive action of our group ${\rm SO}(3)$.

A similar example is discussed in \cite{fels-olver2} that f\/inds dif\/ferential invariants of surfaces in $\R^3$ under the action of ${\rm SO}(3)$, but their invariants are for the {\em variable-parameter} congruence, as opposed to our {\em fixed-parameter} congruence invariants, and hence in particular do not include the f\/irst fundamental form.

\subsection{Another nonhomogeneous space}  Lest we give the wrong impression through the preceding examples that our method is somehow conf\/ined to work only for well-known group actions, now we discuss for our f\/inal example\footnote{Suggested by an anonymous referee.} a~rather randomly chosen intransitive action on $\R^3$ and f\/ind a corresponding complete system of congruence invariants for surface immersions.

Let us consider the following action of a group $G$, parametrized by $(t_1,t_2,t_3,t_4,t_5)$, on $\R^3$ with coordinate system $(u,v,w)$:
\[
G \times \R^3 \longrightarrow \R^3, \qquad (t_1,t_2,t_3,t_4,t_5)\cdot(u,v,w)\longmapsto (u+t_1v+t_2w+t_3,v+t_4w+t_5,w).
\]
(This action can be viewed as one by a certain subgroup of a Heisenberg group as we will see later.)

For the f\/irst step to f\/ind a complete system of congruence invariants for surface immersions $\psi:\R^2\longrightarrow\R^3$ under the action of $G$, we verticalize the given action on $\R^3$ into one on the total space of the product bundle $\R^2\times\R^3\longrightarrow\R^2$ as follows:
\begin{gather*}
 G\times(\R^2\times\R^3)\longrightarrow \R^2\times\R^3,\\
 (t_1,t_2,t_3,t_4,t_5)\cdot(x,y;u,v,w)\longmapsto (x,y;u+t_1v+t_2w+t_3,v+t_4w+t_5,w),
\end{gather*}
where $(x,y)$ represents the coordinate system of $\R^2$, the domain of surface immersions.  This verticalized $G$-action on $\R^2\times\R^3$ is not locally free, but its f\/irst-order prolongation
\begin{gather*}
 \mathcal{A}: \ G\times J^1\big(\R^2,\R^3\big)\longrightarrow  J^1\big(\R^2,\R^3\big),\\
\hphantom{\mathcal{A}: \ } \  (t_1,t_2,t_3,t_4,t_5)\cdot(x,y;u,v,w;u_x,u_y,v_x,v_y,w_x,w_y) \\
\hphantom{\mathcal{A}: \ } \
\longmapsto (x,y;u+t_1v+t_2w+t_3,v+t_4w+t_5,w;\\
\hphantom{\mathcal{A}: \ \longmapsto \ } \ \
   u_x+t_1v_x+t_2w_x,u_y+t_1v_y+t_2w_y,v_x+t_4w_x,v_y+t_4w_y,w_x,w_y),
\end{gather*}
is locally free since we can solve the following chosen normalization equations
\begin{gather*}
 \mathcal{A}^*u=u+t_1v+t_2w+t_3=0, \qquad \mathcal{A}^*v=v+t_4w+t_5=0, \\
  \mathcal{A}^*u_x=u_x+t_1v_x+t_2w_x=0, \qquad \mathcal{A}^*u_y=u_y+t_1v_y+t_2w_y=0, \\
   \mathcal{A}^*v_x=v_x+t_4w_x=0
\end{gather*}
for the group parameters
\begin{gather*}
 t_1=\frac{-u_yw_x+u_xw_y}{v_yw_x-v_xw_y}, \qquad t_2=\frac{u_yv_x-u_xv_y}{v_yw_x-v_xw_y},\\
 t_3=\frac{-uv_yw_x-u_ywv_x+u_yvw_x+uv_xw_y+u_xwv_y-u_xvw_y}{v_yw_x-v_xw_y},\\
 t_4=-\frac{v_x}{w_x}, \qquad t_5=-v+\frac{wz_x}{w_x}.
\end{gather*}
These solutions for $t_1$, $t_2$, $t_3$, $t_4$, $t_5$ are the ones used in the construction of the moving frame
\[
\rho: \ J^1\big(\R^2,\R^3\big)\longrightarrow G, \qquad
(x,y;u,v,w;u_x,u_y,v_x,v_y,w_x,w_y)\longmapsto (t_1,t_2,t_3,t_4,t_5),
\]
and its associated moving frame section $\sigma:J^1(\R^2,\R^3)\longrightarrow G\times J^1(\R^2,\R^3)$.

For the next step of f\/inding Maurer--Cartan forms of $G$, we regard the action of $G$ on $\R^3$ as the linear group consisting of the following elements
\[
g:=\begin{bmatrix}
1 & t_1 & t_2 & t_3 \\
0 & 1 & t_4 & t_5 \\
0 & 0 & 1 & 0 \\
0 & 0 & 0 & 1
\end{bmatrix} \in {\rm GL}(4,\R), \qquad t_i\in \R, \qquad i=1,\dots ,5,
\]
acting on the hyperplane $\{[u \ v \ w \ 1]^{\rm T} \in \R^4 \, | \, u,v,w\in\R\}$ via the usual matrix multiplication.  Then the linearly independent entries of the $\mathfrak{gl(4,\R)}$-valued Maurer--Cartan form $(\dd g)g^{-1}$,  {\samepage
\begin{gather*}
\mu^1:=\dd t_1, \qquad \mu^2:=-t_4\dd t_1+\dd t_2, \qquad \mu^3:=-t_5\dd t_1+\dd t_3, \qquad \mu^4:=\dd t_4, \qquad \mu^5:=\dd t_5,
\end{gather*}
give rise to a Maurer--Cartan coframe f\/ield of $G$.}

Now we use the formula (\ref{constant-structure coframe formula}) to construct constant-structure invariant coframe f\/ield of \linebreak $J^1(\R^2,\R^3)$.  Pulling back the Maurer--Cartan forms of $G$ by the moving frame $\rho$ yields
\begin{gather}
 \rho^*\mu^1=\frac{1}{(v_yw_x-v_xw_y)^2}\Bigl(\! \big(v_yw_xw_y-v_xw_y^2\big)\dd u_x+\big({-}v_yw_x^2+v_xw_xw_y\big)\dd u_y\nonumber\\
\hphantom{\rho^*\mu^1=}{}  +\big({-}u_yw_xw_y+u_xw_y^2\big)\dd v_x
+\big(u_yw_x^2-u_xw_xw_y\big)\dd v_y+(u_yv_xw_y-u_xv_yw_y)\dd w_x\nonumber\\
\hphantom{\rho^*\mu^1=}{}
  +(-u_yv_xw_x+u_xv_yw_x)\dd w_y\!
 \Bigr),\nonumber\\
 \rho^*\mu^2=\frac{1}{v_yw_x^2\!-v_xw_xw_y}\Bigl(\!(-v_yw_x+v_xw_y)\dd u_x+(u_yw_x-u_xw_y)\dd v_x
 +(-u_yv_x+u_xv_y)\dd w_x  \!\Bigr),\nonumber\\
 \rho^*\mu^3=\frac{1}{v_yw_x^2-v_xw_xw_y}\Bigl(\! \big({-}v_yw_x^2+v_xw_xw_y\big)\dd u
+\big(u_yw_x^2-u_xw_xw_y\big)\dd v \nonumber\\
\hphantom{\rho^*\mu^3=}{}
  +(wv_yw_x-wv_xw_y)\dd u_x +(-wu_yw_x+wu_xw_y)\dd v_x
+(-u_yv_xw_x+u_xv_yw_x)\dd w \nonumber\\
\hphantom{\rho^*\mu^3=}{}
  +(wu_yv_x-wu_xv_y)\dd w_x\!\Bigr) ,\nonumber\\
 \rho^*\mu^4=-\frac{1}{w_x}\dd v_x+\frac{v_x}{w_x^2}\dd w_x,\nonumber\\
 \rho^*\mu^5=-\dd v+\frac{v_x}{w_x}\dd w+\frac{w}{w_x}\dd v_x-\frac{wv_x}{w_x^2}\dd w_x,\label{cofMaurerCartan}
\end{gather}
and pulling back the coordinate functions of $J^1(\R^2,\R^3)$, that were not normalized to constants, by $\iota=\mathcal{A}\circ\sigma$ and then taking their exterior dif\/ferentials yields
\begin{gather}
 \dd\iota^*x=\dd x, \qquad \dd\iota^*y=\dd y, \qquad \dd\iota^*w=\dd w, \nonumber\\
 \dd\iota^*v_y=-\frac{w_y}{w_x}\dd v_x+\dd v_y+\frac{v_xw_y}{w_x^2}\dd w_x-\frac{v_x}{w_x}\dd w_y,\nonumber\\
 \dd\iota^*w_x=\dd w_x, \qquad \dd\iota^*w_y=\dd w_y.\label{cofNormalized}
\end{gather}
The eleven one-forms, (\ref{cofMaurerCartan}) and (\ref{cofNormalized}), constitute a constant-structure $G$-invariant coframe f\/ield of $J^1(\R^2,\R^3)$.  The f\/inal step now is to pull back those coframe one-forms by the prolonged graph $j^1\gamma(\psi):\R^2\!\longrightarrow\! J^1(\R^2,\R^3)$ of a generic surface $\psi:\R^2\!\longrightarrow\!\R^3$, $(x,y)\!\longmapsto\! (u(x,y),v(x,y),w(x,y))$, to obtain the following system of one-forms:
\begin{gather}
\big(j^1\gamma(\psi)\big)^*\rho^*\mu^1=\frac{\dd x}{(v_yw_x-v_xw_y)^2}\Bigl(\! -u_y v_{xx} w_x w_y+u_{xx} v_y w_x w_y-u_{xx} v_x w_y^2\nonumber\\
\phantom{\big(j^1\gamma(\psi)\big)^*\rho^*\mu^1=}{} +u_x
   v_{xx} w_y^2 +u_y v_x w_{xx} w_y-u_x v_y w_{xx} w_y-v_y
   w_x^2 u_{x y}+u_y w_x^2 v_{x y}\nonumber\\
\phantom{\big(j^1\gamma(\psi)\big)^*\rho^*\mu^1=}{}
 +v_x w_x w_y
   u_{x y}-u_x w_x w_y v_{x y} -u_y v_x w_x w_{x y}+u_x v_y w_x
   w_{x y} \!\Bigr)\nonumber\\
\phantom{\big(j^1\gamma(\psi)\big)^*\rho^*\mu^1=}{}
  +\frac{\dd y}{(v_yw_x-v_xw_y)^2}\Bigl(\! -u_{yy} v_y w_x^2+u_y v_{yy} w_x^2+u_{yy} v_x w_x w_y-u_x v_{yy} w_x w_y\nonumber\\
\phantom{\big(j^1\gamma(\psi)\big)^*\rho^*\mu^1=}{}
  -u_y v_x w_x w_{yy}+u_x v_y w_x w_{yy}+v_y
   w_x w_y u_{x y}-u_y w_x w_y v_{x y}-v_x w_y^2 u_{x y}\nonumber\\
\phantom{\big(j^1\gamma(\psi)\big)^*\rho^*\mu^1=}{}
 +u_x
   w_y^2 v_{x y}+u_y v_x w_y w_{x y}-u_x v_y w_y w_{x y} \!\Bigr),\nonumber\\
\big(j^1\gamma(\psi)\big)^*\rho^*\mu^2=\frac{\dd x}{v_yw_x^2-v_xw_xw_y}\Bigl(\! u_y v_{xx} w_x-u_{xx} v_y w_x-u_y v_x w_{xx}+u_x v_y
   w_{xx}\nonumber\\
\hphantom{\big(j^1\gamma(\psi)\big)^*\rho^*\mu^2=}{}
  +u_{xx} v_x w_y-u_x v_{xx} w_y \Bigr)+\frac{\dd y}{v_yw_x^2-v_xw_xw_y}\Bigl(\! -v_y w_x u_{x y}+u_y w_x v_{x y}\nonumber\\
\hphantom{\big(j^1\gamma(\psi)\big)^*\rho^*\mu^2=}{}
+v_x w_y u_{x   y}  -u_x w_y v_{x y}-u_y v_x w_{x y}+u_x v_y w_{x y} \!\Bigr),
\nonumber\\
\big(j^1\gamma(\psi)\big)^*\rho^*\mu^3=\frac{\dd x}{v_yw_x^2-v_xw_xw_y}\Bigl(\! -w u_y v_{xx} w_x+w u_{xx} v_y w_x+w u_y v_x w_{xx}\nonumber\\
\hphantom{\big(j^1\gamma(\psi)\big)^*\rho^*\mu^3=}{}
 -w u_x v_y w_{xx}-w u_{xx} v_x
   w_y+w u_x v_{xx} w_y \!\Bigr)+\frac{\dd y}{v_yw_x^2-v_xw_xw_y}\Bigl(\! w v_y w_x u_{x y}\nonumber\\
\hphantom{\big(j^1\gamma(\psi)\big)^*\rho^*\mu^3=}{}
  -w u_y w_x v_{x y}-w v_x w_y u_{x y}+w u_x w_y v_{x y}+w u_y v_x w_{x
   y}-w u_x v_y w_{x y}\! \Bigr), \nonumber\\
\big(j^1\gamma(\psi)\big)^*\rho^*\mu^4=\Bigl(\!\frac{v_x w_{xx}}{w_x^2}-\frac{v_{xx}}{w_x}\!\Bigr)\dd x+\Bigl(\!\frac{v_x w_{x y}}{w_x^2}-\frac{v_{x y}}{w_x}\!\Bigr)\dd y, \nonumber\\
\big(j^1\gamma(\psi)\big)^*\rho^*\mu^5=\Bigl( \!\frac{w v_{xx}}{w_x}-\frac{w v_x w_{xx}}{w_x^2} \!\Bigr)\dd x+\Bigl(\! \frac{v_x w_y}{w_x}-\frac{w v_x w_{x y}}{w_x^2}+\frac{w v_{x y}}{w_x}-v_y \!\Bigr)\dd y, \nonumber\\
\big(j^1\gamma(\psi)\big)^*\dd\iota^*x=\dd x, \nonumber\\
\big(j^1\gamma(\psi)\big)^*\dd\iota^*y=\dd y, \nonumber\\
\big(j^1\gamma(\psi)\big)^*\dd\iota^*w=w_x\dd x+w_y\dd y, \nonumber\\
\big(j^1\gamma(\psi)\big)^*\dd\iota^*v_y=\Bigl(\! \frac{v_x w_{xx} w_y}{w_x^2}-\frac{v_{xx} w_y}{w_x}-\frac{v_x w_{x y}}{w_x}+v_{x y} \!\Bigr)\dd x+\Bigl( \! -\frac{v_x w_{yy}}{w_x}\nonumber\\
\hphantom{\big(j^1\gamma(\psi)\big)^*\dd\iota^*v_y=}{}
 -\frac{w_y v_{x y}}{w_x}+\frac{v_x w_y w_{x y}}{w_x^2}+v_{yy}\! \Bigr)\dd y, \nonumber\\
\big(j^1\gamma(\psi)\big)^*\dd\iota^*w_x=w_{xx}\dd x+w_{xy}\dd y, \nonumber\\
\big(j^1\gamma(\psi)\big)^*\dd\iota^*w_y=w_{xy}\dd x+w_{yy}\dd y.\label{final pullback of G-invariant coframe}
\end{gather}
Thus, a complete system of $G$-congruence invariants for surface immersions in $\R^3$ is provided by the set of all the coef\/f\/icients of $\dd x$ and $\dd y$ in (\ref{final pullback of G-invariant coframe}).  According to Theorem~\ref{general solution}, two surface immersions will be congruent under the action of $G$ if and only if some transformation $g\in G$ takes some f\/irst-order jet of one of the surfaces to a jet of the other surface and the inva\-riants~(\ref{final pullback of G-invariant coframe}) evaluated for the two surfaces agree.

\section{Discussions}

We introduced the notion of invariant coframe f\/ields of constant structure, whose explicit construction can be done by the method of equivariant moving frames, and used them to prove  Theorem~\ref{general solution} that provided theoretical justif\/ication of our coherent method of completely sol\-ving the congruence problem of immersions in homogeneous and nonhomogeneous spaces alike.
It extends and generalizes to arbitrary $G$-spaces the key congruence lemmas in \cite{griffiths} that were designed for homogeneous spaces.
We demonstrated our method by applying it to congruence problems in some classical and other examples.

The next order of research in this direction should involve applications of our method to more substantial and unexplored congruence problems.

\subsection*{Acknowledgments}
This work has benef\/ited from the discussions held in the Dif\/ferential Geometry and Lie Theory seminars  at the University of Toledo; the author would like to thank the organizers and parti\-ci\-pants of the seminars.  Also, the anonymous referees' critical and yet helpful comments have contributed signif\/icantly in the process of revising and improving the paper;  the author is very grateful to the referees.

It is hoped that this work serves to ref\/lect, although only to a small extent limited by the \mbox{author's} meager knowledge, the author's appreciation of the introduction by Professor Peter Olver to the marvelous unifying philosophy and technology of symmetry, invariance, and equiva\-lence.

\pdfbookmark[1]{References}{ref}
\LastPageEnding

\end{document}